\documentclass[11pt,final,a4paper]{amsart}
\usepackage[latin1]{inputenc}
\usepackage[english]{babel}
\usepackage{amsmath}
\usepackage{amsfonts}
\usepackage{amssymb}
\usepackage{graphicx}
\usepackage{indentfirst}
\usepackage{dsfont}
\usepackage{color}
\usepackage{enumitem}

\newtheorem{theorem}{Theorem}[section]

\newtheorem{proposition}[theorem]{Proposition}
\newtheorem{corollary}[theorem]{Corollary}
\title{Hypothesis testing for Markovian models with random time observations}
\date{}

\newcommand{\vect}{\operatorname{vec}}
\newcommand{\id}{\operatorname{I}}
\newcommand{\supp}{\operatorname{supp}}

\newcommand{\lin}{\operatorname{lin}}
\newcommand{\Com}{\operatorname{Com}}
\newcommand{\im}{\operatorname{Im}}

\newcommand{\rank}{\operatorname{rank}}
\newcommand{\Vvert}{\Vert \hspace{-0.04cm} \vert}

\begin{document}
\maketitle 

\begin{center}
Flavia Barsotti\footnote{Risk Methodologies, Group Financial Risks, Group Risk Management, UniCredit S.p.A, 20154 Milano\\
The views expressed in this paper are those of the author and should not be attributed to UniCredit Group or to the author as representative or employee of UniCredit Group.}, Anne Philippe\footnote{Universit\'e de Nantes, Laboratoire de math\'ematiques Jean Leray} and Paul Rochet{\footnotesize${}^{ 2}$}
\end{center}

\begin{abstract}
The aim of this paper is to propose a methodology for testing general hypothesis in a Markovian setting with random sampling. A discrete Markov chain $X$ is observed at random time intervals $\tau_k$, assumed to be iid with unknown distribution $\mu$. Two test procedures are investigated. The first one is devoted to testing if the transition matrix $P$ of the Markov chain $X$ satisfies specific affine constraints, covering a wide range of situations such as symmetry or sparsity. The second procedure is a goodness-of-fit test on the distribution $\mu$, which reveals to be consistent under mild assumptions even though the time gaps are not observed. The theoretical results are supported by a Monte Carlo simulation study to show the performance and robustness of the proposed methodologies on specific numerical examples. 
\end{abstract}

\vspace{3mm}


\normalsize

\section{Introduction}
\label{sec:int}

It has been widely recognized that discrete Markov chains are a powerful probabilistic tool to study many real phenomena in different application fields. Interests from economics, finance, insurance and also medical research are just few examples. From a theoretical point of view, alternative technical frameworks have been defined and developed in the literature to perform statistical inference in Markovian models: multiple Markov chains \cite{MR0084903,MR0039207}, hidden Markov processes \cite{MR0202264,cappe2009inference,rabiner1989tutorial}, random walks on graphs \cite{gkantsidis2004random} are well known mathematical settings used to describe the evolution of real events. \\

In this paper, we investigate hypothesis testing issues in a Markov model with random sampling. In the spirit of \cite{fytp14}, a discrete homogenous Markov chain $(X_m)_{m \in \mathbb N}$ is observed at random times so that the only available observations consist in a sub-sequence of the initial process. Hereafter we denote by $(Y_k)_{k \in \mathbb N} $ the observed process.  The time gaps $\tau_k$ (i.e. the number of jumps) between two consecutive observations are assumed non-negative, independent and identically distributed from an unknown distribution $\mu$. A statistical methodology is proposed to test specific hypothesis on the transition matrix $P$ or the distribution $\mu$ having observed $(Y_1,...,Y_n)$. \\

Our interest is to answer the following questions: can we test some specific hypothesis on the transition matrix $P$ of the initial chain $X$ when neither the time gaps $\tau_k$ nor their distribution $\mu$ are known? Additionally, can we perform a goodness-of-fit test on the distribution $\mu$ of the time gaps in this framework?\\

As explained in \cite{fytp14}, a typical application of this setting occurs when considering a continuous time Markov chain $Z=(Z_t)_{t>0}$ observed at a discrete time grid $t_1 < ... < t_n$.  In this situation, $X$ represents the jump process of $Z$  and $(\tau_k)_{k \in \mathbb N}$ denotes the number of jumps occurring between two consecutive observations, $Y_k:=Z_{t_k}$. If the discrete time grid is chosen independently from the chain, the time gaps $\tau_k$ are independent random variables {\sl{unknown}} to the practitioner. Their distribution $\mu$ is Poisson if $Z$ is observed at  regularly spaced time intervals, but one can easily imagine a different distribution $\mu$ if $t_1, ...t_n$ are affected by undesired random effects. This framework can describe many real phenomena: chemical reactions \cite{anderson2011continuous}, financial markets \cite{MR1822778}, waiting lines in queuing theory \cite{gaver1959imbedded}, medical studies \cite{craig02}.\\

In \cite{fytp14}, the sparsity of $P$ reveals to be a crucial assumption for the proposed estimation methodology. In this paper, we extend the statistical setting defined in \cite{fytp14} by relaxing the sparsity assumption and working under a more general framework where the information about $P$ can be expressed in the form of affine constraints. It is then possible to incorporate some prior information on $P$ such as symmetry, reflexivity or even simply the knowledge of some of its entries. The sparsity hypothesis considered in \cite{fytp14} can be retrieved as a particular case. \\

In our model, the observations $Y_1,...,Y_n$ are the realizations of an homogenous Markov chain, with a transition probability $Q$ that can be written as an analytic function of $P$. As in \cite{fytp14}, neither the transition matrix $P$ nor the distribution $\mu$ are known and the time gaps $\tau_k$ between two consecutive observations $Y_{k-1}, Y_k$ are unavailable. Our aim is now to perform hypothesis tests on both the transition kernel $P$ and the time gaps' distribution $\mu$. Theoretical results are separated in two consecutive steps. As a first step, we describe how to build an hypothesis test on $P$ to study if the transition kernel $P$ satisfies additional affine constraints. This framework can be used to test a wide range of hypothesis such as symmetry,  sparsity or even to test a particular value of $P$ or some of its entries. As a second step, we propose a procedure to test specific values for the distribution of time gaps $\mu$, considering an hypothesis of the form $\mathcal H_0: \mu = \mu_0$, given some preliminary information on $P$. Both tests are asymptotically exact and rely on the asymptotic normality of the empirical transition matrix $\widehat Q$. A Monte Carlo simulation analysis is performed in order to highlight the performance of the proposed test.\\

The paper is organized as follows: Section \ref{sec:prob} provides a detailed description of the statistical problem and discusses the identifiability issues encountered in this framework. Section \ref{sec:H} defines an asymptotically exact hypothesis test for the transition kernel $P$. Section \ref{sec:mu} shows the theoretical framework to build an hypothesis test on the distribution $\mu$. Section \ref{sec:cs} supports the analysis of Sections \ref{sec:H} and \ref{sec:mu} with numerical results from a Monte Carlo simulation study. Proofs are postponed in the Appendix.

\section{The problem}\label{sec:prob}

We consider an irreducible homogenous Markov chain $X=(X_m)_{m \in \mathbb N}$ with finite state space $\mathcal E = \{1,...,N\}$, $N \geq 3$ and transition matrix $P$. We assume that $X$ is observed at random times $T_1,...,T_n$ so that the only available observations consist in the sub-sequence $Y_k:= X_{T_k}, k=1,...,n$. The numbers of jumps $\tau_k:= T_{k-1} - T_k$ between two observations $Y_k$ are assumed to be iid random variables with distribution $\mu$ on $\mathbb N$ and independent from $X$. In this setting, the resulting process $Y=(Y_k)_{k \in \mathbb N}$ remains Markovian (see \cite{fytp14})
with transition matrix
\begin{equation}
Q:= G_\mu(P) = \sum_{\ell\geq 0} P^\ell \mu(\ell),\label{ap:2}
\end{equation}
where $G_\mu:[-1,1] \to \mathbb R$ denotes the generator function of $\mu$. We assume that some information is available on $P$ in the form of an affine constraint
$$ A p = b    $$
where $p = \vect(P) = (P_{11},...,P_{N1},...,P_{1N},..., P_{NN})^\top$ is the vectorization of $P$, $A \in \mathbb R^{r \times N^2}$ is a known full ranked matrix and $b \in \mathbb R^{r}$ a known vector. This affine constraint defines the model 
\begin{equation}
\mathcal M := \{ m  \in \mathbb R^{N^2}: A m = b \},  \label{ap:1}
\end{equation}
that is, the set of admissible values for $p=\vect(P)$. Each additional affine condition satisfied by $P$ reduces by one the dimension of the model. Therefore, the rank $r$ of $A$ indicates the information available on $P$. In the sequel, we denote by $d$ the dimension of the model $\mathcal M$, given by $ d =\dim(\mathcal M) = N^2 - r$.\\

Many different forms of information on the initial chain $X_m$ can be expressed by affine constraints, as described in the following examples. 
\begin{enumerate}[label=\alph*)]
	\item $ P \mathbf 1 = \mathbf 1$, where $\mathbf 1 = (1,...,1)^\top$, is an important information that can be expressed as an affine constraint on $P$ by
$$ (\mathbf 1^\top \otimes \id ) \vect(P)  = \mathbf 1, $$
where $\otimes$ stands for the Kronecker product.
  \item In \cite{craig02}, the author investigates the progression of a disease using a Markov chain observed at unequal time intervals. For this particular application, the transition matrix of the chain is expected to be triangular: in this case the information can be represented with $N(N-1)/2$ linear conditions on $P$.
	\item If the initial chain $X$ cannot remain at the same state on two consecutive times, the transition matrix $P$ is known to have zero diagonal. This information can be expressed in the form of $N$ linear constraints on $P$.
  \item More generally, if one or several state transitions are known to be impossible in the initial chain $(X_m)_{m \in \mathbb N}$, the nullity of the corresponding entries can be expressed as a set of linear conditions on $p$. The number of constraints is equal to the number of known zeros in $P$. This situation is treated in \cite{fytp14}.
	\item A symmetric transition matrix $P$ corresponds to $N(N-1)/2$ linear constraints on $p = \vect(P)$. 
	\item $P$ is doubly stochastic, one can include the condition $P^\top \mathbf 1 = \mathbf 1$ as $N-1$ additional affine constraints. 
	\item $X$ is a reversible Markov chain, its transition matrix $P$ satisfies the linear conditions $\pi_i P_{ij} = \pi_j P_{ji}$ for all $i,j=1,...,N$, where $\pi=(\pi_1,...,\pi_N)$ is the invariant measure. Although $\pi$ is unknown in practice, it can be estimated directly from the observations since it is also the invariant measure of $Q$. 
	\item The knowledge of some entries of $P$ being equal, or equal to a known value 
	$ c \in [0,1]$ can be expected in certain practical situations (e.g. the transition probability distribution from a given state $i$ is uniform). These are of course particular examples of affine constraints on $P$.
\end{enumerate}

The positivity of the entries of $P$ appears to be more difficult to handle than affine constraints. Nevertheless, the positivity is somewhat less informative as it does not reduce the dimension of the model. In this paper, we choose to neglect this information on $P$ for simplicity. \\

To avoid critical situations, we assume that $(X_m)_{m \in \mathbb N}$ is an aperiodic Markov chain, in which case the transition kernel $P$ has a unique invariant distribution $\pi = (\pi_1,...,\pi_N)$, where $\pi_i$ is positive for all $i=1,...,N$. From the relation $Q = G_\mu(P)$, we know that, like $P$, $Q$ is aperiodic recurrent and they share the same invariant distribution. Thus, $\pi$ can be estimated consistently from the observations $Y_k$, along with $Q$ by the empirical estimators
\begin{equation}\label{eq:Qchap}  \hat \pi_i = \frac 1 n \sum_{k=1}^n \mathds 1 \{ Y_k = i \} \ \ , \ \ \hat Q_{ij} = \frac{\sum_{k=1}^{n-1} \mathds 1 \{ Y_k=i,Y_{k+1}=j\}}{ \sum_{k=1}^{n-1} \mathds 1 \{ Y_k=i\} }, \end{equation}
for all $i,j =1,...,N$, provided that the state $i$ has been observed at least once. It is known (see for instance \cite{guttorp1995stochastic}) that $\hat Q$ is  asymptotically unbiased  and  $\hat q =\vect(\hat Q) $
is  asymptotically Gaussian,  
$$ \sqrt n (\hat q - q) \stackrel{\ d \ }{\longrightarrow}  \mathcal N(0, \Sigma). $$ 
The matrix $\Sigma  \in \mathbb R^{N^2 \times N^2}$ can be deduced from the asymptotic behavior of the covariances 
\begin{equation}\label{empi} \forall i,j,k,\ell = 1,...,N, \   \lim_{n \to \infty} n \operatorname{cov}(\hat Q_{ij}, \hat Q_{k\ell})  = \left\{ \begin{array}{cl} Q_{ij}(1 - Q_{ij})/\pi_i & \operatorname{if } \ (i,j) = (k,\ell), \\
 - Q_{ij} Q_{i\ell}/\pi_i & \operatorname{if } \ i=k, \ j\neq \ell, \\
0 & \operatorname{otherwise.} \end{array} \right. \end{equation}
Remark that $\Sigma$ is a singular matrix since $\hat q - q$ can only take values in a linear subspace of $\mathbb R^{N^2 \times N^2}$. 

Since the distribution of the time gaps $\tau_k$ is unknown, the main information available on $P$ is that it commutes with $Q$.
Thus, given the model $\mathcal M$, the possible values for $p$ can be restricted to $\mathcal M \cap \Com(Q)$, where $\Com(Q)$ is the commutant of $Q$. As in \cite{fytp14}, we consider the commutation operator $m \mapsto \Delta(Q) m = \vect(MQ - QM) $ where 
$$ \Delta(Q) =  \id \otimes \ Q - Q^\top \otimes \id.$$  
Of course, since $Q$ is unknown, we use the empirical estimator $\hat Q$ to retrieve some information on $P$.

\section{Hypothesis tests on $P$}\label{sec:H}

Given that the transition matrix $P$ of the initial Markov process $(X_m)_{m \in \mathbb N}$ belongs to some model $\mathcal M$ of the form $\eqref{ap:1}$, we want to test additional affine conditions on $p= \vect(P)$, via the null hypothesis $\mathcal H_0: A_0 p  = b_0$, where $A_0 \in \mathbb R^{k \times N^2}$ is a full ranked matrix, $b_0 \in \mathbb R^{k}$ and $k \geq 1$ is the number of new conditions to test. To avoid technical issues, the new constraints must be compatible with the model, i.e. there exits at least one element $m\in \mathcal M$ such that $A_0 m = b_0$. Moreover, we assume without loss of generality that the new constraints are not redundant, in the sense that each row of $A_0$ brings a new information on $P$. Formally, this means that the model under $\mathcal H_0$,
$$ \mathcal M_0 := \{ m \in \mathbb R^{N^2}: A m = b, A_0 m = b_0 \} \subset \mathcal M, $$
is of dimension $d - k$ where we recall $d = \dim(\mathcal M)$. Let $\phi_1,...,\phi_{d-k}$ be a basis of $\ker(A) \cap \ker(A_0)$, extended by $\phi_{d-k+1},...,\phi_d$ to form a basis of $\ker(A)$. The initial model $\mathcal M$ and the constrained model $\mathcal M_0$ can be expressed as
$$\mathcal M = \{ p_0 + \Phi \beta: \beta \in \mathbb R^d \} \ \ \text{ and } \ \ \mathcal M_0 = \{ p_0 + \Phi_0 \gamma: \gamma \in \mathbb R^{d-k} \},$$ 
setting $\Phi = ( \phi_1,...,\phi_d) \in \mathbb R^{N^2 \times d}$, $\Phi_0 =  (\phi_1,...,\phi_{d-k}) \in \mathbb R^{N^2 \times (d-k)}$ and with $p_0$ an element of $\mathcal M_0$. Define the linear sets
\begin{eqnarray*} E & = & \im(\Delta(Q) \Phi) = \{\Delta(Q) m: m \in \ker(A) \} \\
 E_0 \!\! & = & \im(\Delta(Q) \Phi_0) = \{\Delta(Q) m: m \in \ker(A) \cap \ker(A_0) \} \\ 
F & = & E \cap E_0^\perp,
\end{eqnarray*}
where $E_0^\perp$ stands for the orthogonal complement of $E_0$ in $\mathbb R^{N^2}$. Note that since $E_0 \subseteq E$, we have $E = E_0 \oplus F$. Moreover, we shall denote by $\Pi_E$, $\Pi_{E_0}$ and $\Pi_F$ the orthogonal projectors in $\mathbb R^{N^2}$ onto $E, E_0$ and $F$ respectively. We make the following assumptions.
\begin{itemize}
	\item[\textbf{A1.}] The problem is identifiable under $\mathcal H_0$, i.e. $\Com(Q) \cap \ker(A) \cap \ker(A_0) = \{ 0 \}$.
	\item[\textbf{A2.}] The dimension of $E$ is maximal over the set of stochastic matrices with the same support as $Q$:
	$$ \dim(E) = \rank(\Delta(Q) \Phi ) = \max_{\substack{Q': Q' \mathbf 1 = \mathbf 1 \\ \supp(Q') \subseteq \supp(Q)}} \rank(\Delta(Q') \Phi).   $$
\end{itemize}

The equivalence, stated in \textbf{A1}, between the identifiability and the matrix $\Delta(Q) \Phi_0$ being of full rank is mentioned in \cite{fytp14}. Here, the identifiability under $\mathcal H_0$ means that $P$ is fully characterized by the condition $ p \in \mathcal M_0$ and the fact that $P$ and $Q$ commute. Theorem 3.2 in \cite{fytp14} ensures the existence of a consistent estimator of $p$ as soon as the hypothesis $\mathcal H_0$ is true. Remark that unlike in \cite{fytp14}, we do not require the initial model $\mathcal M$ to be identifiable, although if it is the case, the two assumptions always hold. 

Though the assumption {\textbf{A2} is  quite technical, it is very mild, since the set of matrices $Q'$ for which the rank of $\Delta(Q') \Phi$ is maximal is a dense open set.} The argument is similar to the one used in the proof of Lemma 2.1 in \cite{fytp14}.  \\

Hereafter we denote by $\chi^2_G(W)$ the centered generalized chi-squared distribution
$$ \chi^2_G(W) \overset{d}{=} \epsilon^\top W \epsilon,   $$
where $\epsilon$ is a standard Gaussian vector in $\mathbb R^{N^2}$ and $W\in \mathbb R^{N^2\times N^2} $ . Remark that the above distribution is simply a Dirac mass at zero if $W$ is the null matrix. 

\begin{theorem}\label{th1} Suppose that \textbf{A1} and \textbf{A2} hold. Under $\mathcal H_0: A_0 p  = b_0$, as $n\to \infty$ the statistic
$$ S = n \bigg( \inf_{m \in \mathcal M_0} \Vert \Delta(\hat Q) m \Vert^2 - \inf_{m \in \mathcal M} \Vert \Delta(\hat Q) m \Vert^2 \bigg)  $$
converges in distribution to a $\chi^2_G(W)$ as $n \to \infty$ where 
\begin{equation}
W : = \Pi_F \Delta(P) \Sigma\Delta(P)^\top\Pi_F\label{ap:3}
\end{equation}
\end{theorem}
\noindent The proof is given in Appendix. \\ 

This result does not rule out the possibility that $W$ is the null matrix (the limit distribution is $0$ in this case), however, it is of no use to perform the test in this situation.\\

Theorem \ref{th1} alone is not sufficient to build a test for $\mathcal H_0$ since the limit distribution is unknown in practice. Nevertheless, $W$ can be estimated consistently under $\mathcal H_0$ using estimates of $\Pi_F$, $P$ and $\Sigma$. 
\begin{itemize}[label=-]
	\item An estimator of $\Sigma$ can be easily obtained by using the plug in principle. We replace in \eqref{empi} the parameters $\pi$ and $Q$ by the empirical versions $\hat \pi$ and $\hat Q$ defined in \eqref{eq:Qchap} . Since $\pi_i$ is positive for all $i=1,...,N$, the resulting estimator $\hat \Sigma$ is consistent. 
	\item In \cite{fytp14}, the authors show that a consistent estimator of $P$ can be obtained if the problem is identifiable. In our framework, the identifiability under $\mathcal H_0$ assumed in \textbf{A1} reduces to
	$$\{ P \} = \mathcal M_0 \cap \Com(Q). $$	
From Lemma 6.1 in \cite{fytp14}, this condition is actually equivalent to $\Delta(Q) \Phi_0$ being of full rank, in which case
\begin{equation}\label{eq:pchap} \hat p =  \left( \id - \Phi_0 \big(\Phi_0^\top \Delta(\hat Q)^\top \Delta(\hat Q) \Phi_0 \big)^{-1} \ \Phi_0^\top \Delta(\hat Q)^\top \Delta(\hat Q) \right) p_0  \end{equation}
is a consistent estimator of $p$ under $\mathcal H_0$.
\item In the proof of Theorem \ref{th1} (see Appendix)  an estimator of $\Pi_F$ is  obtained by considering the linear space 
$$ \hat F = \im(\Delta(\hat Q)\Phi) \cap \im(\Delta(\hat Q)\Phi_0)^\perp.    $$
Under \textbf{A1} and \textbf{A2}, we prove that  the resulting projector $\Pi_{\hat F}$ converges in probability to $\Pi_F$ as $n \to \infty$.
\end{itemize}
By plugging the  estimates $\hat \Sigma$, $\hat P$ and $\Pi_{\hat F}$ in   \eqref{ap:3}, we  build 
$$ \hat W :=  \Pi_{\hat F} \Delta(\hat P ) \hat \Sigma \Delta(\hat P )^\top \Pi_{\hat F}. $$
which is a consistent estimate of $W$ under $\mathcal H_0$. As a result, the limit distribution of the statistic $S$ can be approximated by $\chi^2_G(\hat W)$. This leads to the following result.

\begin{corollary}\label{coro} Assume that \textbf{A1} and \textbf{A2} hold. Let $\hat u_{1-\alpha}$ denote the quantile of order $1- \alpha$ of the $\chi^2_G(\hat W)$ distribution for $\alpha \in (0,1)$. If $W\neq 0$, we have under $\mathcal H_0: A_0 p = b_0$, 
$$ \lim_{n \to \infty} \mathbb P(S > \hat u_{1-\alpha}) = \alpha. $$
Moreover, if $\forall p' \in \mathcal M_0: \ \Delta(Q) p' \neq 0$, then $\displaystyle \lim_{n \to \infty} \mathbb P(S > \hat u_{1-\alpha}) = 1$.
\end{corollary}
\noindent The proof is given in Appendix.\\

Because the quantiles of the generalized $\chi^2$ distribution are not available in an analytic form,  they are approximated by Monte-Carlo methods, following \cite{davies1980algorithm}.

\section{Test on $\mu$}\label{sec:mu}
The estimation of the transition kernel $P$ is based on the property that $P$ and $Q$ commute, which is due to the existence of an analytic function $ G_\mu(\cdot)$ such that $Q = G_\mu(P)$ defined in \eqref{ap:2}. However, no particular interest has been given so far to the actual value of $G_\mu(\cdot)$, since the knowledge of the distribution $\mu$ is not needed to build the test procedure on $\hat P$.\\

We emphasize that in our framework it is in general not possible to exactly recover the distribution $\mu$, since different distributions $\mu$ and $\nu$ may lead to the same image $G_\nu(P)=G_\mu(P)$. Nevertheless, questions regarding the number of jumps between consecutive observations, or their distribution $\mu$ can naturally arise in practical cases. Suppose for instance that the observations $Y_i$ come from a continuous process $(X_t)_{t \geq 0}$ observed at regular time intervals. Assuming that nothing other than the $Y_i$'s is known, one might be interested in checking if the underlying process $(X_t)_{t \geq 0}$ is a continuous time Markov chain. This hypothesis relies on the distribution of the time intervals between two consecutive jumps of $X_t$. If the process $(X_t)_{t \geq 0}$ is in fact Markovian, the times between jumps are drawn from an exponential distribution, resulting in Poisson random variables $\tau_k$. In such a case, testing this hypothesis on $\mu$ reduces to verify if $G_\mu$ is the generator function of a Poisson variable, that is, if $G_\mu(t) = e^{\lambda(t-1)}$ for some $\lambda >0$.\\

In this section we  consider  hypotheses of the form $\mathcal H_0: \mu = \mu_0$. We assume that the model is identifiable, in which case we can build a consistent estimator of $p$ as follows 
$$  \hat p =  \left( \id - \Phi \big(\Phi^\top \Delta(\hat Q)^\top \Delta(\hat Q) \Phi \big)^{-1} \ \Phi^\top \Delta(\hat Q)^\top \Delta(\hat Q) \right) p_0.$$ 
From Theorem 3.2 in \cite{fytp14}, we know that $\hat p$ satisfies
\begin{equation}\label{eqhatp} \sqrt n \ \big( \hat p - p \big) = \sqrt n \ B \big( \hat q - q \big) + o_P ( 1 ), \end{equation}
where $B = \Phi \left[ \Phi^\top  \Delta(Q)^\top \Delta(Q) \Phi \right]^{-1} \Phi^\top \Delta(Q)^\top \Delta(P)$.\\
 
The procedure to perform a test on the hypothesis $\mathcal H_0: \mu = \mu_0$ relies on the comparison between $\hat Q$ and $G_{\mu}(\hat P)$. Since $ Q=G_{\mu}(P)$, the consistency of $\hat P$ implies that $\hat Q  - G_{\mu}(\hat P)$ converges in probability to zero. Moreover, we can derive the asymptotic distribution of $g_{\mu}(\hat P) := \vect(G_{\mu}(\hat P))$, leading to the following result. 

\begin{proposition}\label{test} If the time gaps distribution $\mu$ satisfies the moment condition $\mathbb E_\mu(\tau) < \infty$, then the matrix $ \Gamma =  \sum_{k \geq 1} \big( \textstyle \sum_{j=1}^{k} (P^{j-1})^\top \otimes P^{k-j}\big) \mu(k)$ is well defined and 
$$ \sqrt n \big(\hat q - g_{\mu}(\hat P) \big) \stackrel{\ d \ }{\longrightarrow}  \mathcal N \Big(0, \big(\id - \Gamma B \big) \Sigma \big( \id - \Gamma B \big)^\top \Big).   $$
\end{proposition}
\noindent The proof is given in Appendix. \\
Consequently we get,  under $\mathcal H_0$,  
$$ n ||\hat q - g_{\mu_0}(\hat P)  ||^2 \stackrel{\ d \ }{\longrightarrow}   \chi^2_G\big( (\id - \Gamma  B )  \Sigma ( \id - \Gamma B)^\top \big).$$ 
This result  is not sufficient to build the statistical test on the distribution $\mu$. So, we proceed similarly as before by using consistent estimates of $\Gamma$, $B$ and $\Sigma$ to replace the unavailable true values in order to approximate the limit variance. Then, the resulting test statistic has 
critical region of the form $\{n ||\hat q - g_{\mu_0}(\hat P)  ||^2 > \tilde u_{1-\alpha} \}$ where $\tilde u_{1-\alpha}$  is the quantile of the $\chi^2_G\big( (\id - \hat\Gamma \hat B ) \hat\Sigma ( \id - \hat\Gamma \hat B)^\top \big)$ distribution. 
 
\section{Computational study}
\label{sec:cs}
This section is devoted to a Monte Carlo simulation analysis for three hypothesis tests proposed in this paper.  In the simulations, the initial Markov chain $X$ behaves like a reflected random walk. The chain takes values in some finite space having $N=10$ states. Our ultimate aim  is  to test the  assumption $P = P_0$, where $P_0$ is  given by
\begin{equation}
\label{ex:p1} P_0 = \left[ \begin{array}{ccccc} 
0 \! & \! 1 \! & \! 0 \! & \! \dots \! & \!  \\
0.5 \! & \! 0 \! & \! 0.5 \! & \! \ddots \! & \! \vdots \\
0 \! & \! \ddots \! & \! \ddots \! & \! \ddots \! & \! 0 \\
\vdots \! & \! \ddots  \! & \! 0.5 \! & \! 0 \! & \! 0.5 \\
 \! & \! \dots \! & \! 0 \! & \! 1 \! & \! 0 \\
 \end{array} \right].  
 \end{equation}

Assuming that nothing is known on the transition kernel $P$, we divide the hypothesis tests on $P$ into three steps. First, we want to test if the support of $P$ is contained in that of $P_0$, i.e. if only the upper and lower diagonals of $P$ have non-zero entries. Secondly, we investigate the test procedure for the null hypothesis $P=P_0$, under the assumption that the support of $P$ is restricted to that of $P_0$. Finally, a third test is carried out for hypothesis on the distribution $\mu$. \\

We perform the three tests for four sample sizes: $n=200$, $n=500$, $n=1000$ and $n=2000$. Under the null hypothesis, we simulate the Markov chain $X=(X_m)_{m \in \mathbb N}$ with transition matrix $P_0$ and the times gaps $\tau_k$ are drawn from a Poisson distribution $\mu = \mathcal P(1)$. Thus, the available observations are $Y_k:= X_{T_k}$ where $T_k = \sum_{j=1}^k \tau_k$ for $k=1,\dots n$. We recall that neither the initial Markov chain $X$ nor the time gaps $\tau_k$ are observed and the distribution $\mu$ of the $\tau_k$ is unknown to the practitioner. The nominal significance level is fixed equal to $\alpha =5\%$ for all numerical experiments.

\subsection*{Test 1.}
We wish to test if the support of $P$ is contained in that of $P_0$, i.e. if only the upper and lower diagonals of $P$ have non-zero entries.
The hypothesis  is
$$ \mathcal H_0^{(1)}: \supp(P) \subseteq \supp(P_0),  $$
and the model $\mathcal M^{(1)}$ is maximal,  i.e. it contains all stochastic matrices: $\mathcal M^{(1)} = \{ \vect(M): M \mathbf 1 = \mathbf 1 \}$. By Monte-Carlo method, the significance level is estimated for the four sample sizes using  $10^4$ replications.  Results are gathered in Table \ref{tab:exemple1}. 

\begin{table}[ht]
\begin{center}
\begin{tabular}{|c||c|c|c|c|} \hline 
$ n $ & $200$ & $500$ & $1000$ & $2000$ \\
\hline
\hline
$\hat \alpha$ & $\!\!\! \begin{array}{c} 0.131 \\ (0.0034) \end{array}\!\!\! $ & $ \!\!\! \begin{array}{c} 0.064 \\ (0.0024) \end{array}\!\!\!  $ & $\!\!\! \begin{array}{c} 0.059 \\ (0.0024) \end{array}\!\!\! $ & $\!\!\! \begin{array}{c} 0.056 \\ (0.0023) \end{array}\!\!\! $ \\
\hline
\end{tabular}
\vspace{0.3cm}
\caption{\footnotesize\textit{Test  for hypothesis  $\mathcal H_0^{(1)}: \supp(P) \subseteq \supp(P_0)$ in the model $\mathcal M^{(1)} = \{ \vect(M): M \mathbf 1 = \mathbf 1 \}$.} The table contains summary statistics of Monte Carlo simulation results obtained with $P=P_0$ given in Eq. (\ref{ex:p1}), for the four different sample sizes $n$ and for $\mu = \mathcal P(1)$. The probability to reject $\mathcal H_0^{(1)}$ is estimated for a nominal size  $\alpha=0.05$, based on $10^4$ replications. Standard errors are given in brackets.} 
\label{tab:exemple1}
\end{center}
\end{table} 

\vspace{-0.6cm}

\noindent Table \ref{tab:exemple1} shows that the estimated significance level converges to the nominal size $5\%$  as $n$ increases. It seems that the test is  accurate as soon as the sample sizes is  greater than 500.\\

\begin{figure*}[ht]
		\includegraphics[width=0.46\textwidth,height=4.5cm]{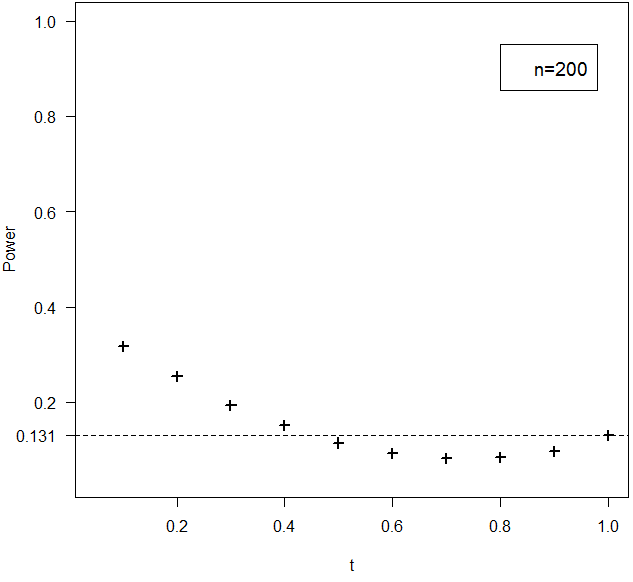} \hfill \includegraphics[width=0.46\textwidth,height=4.5cm]{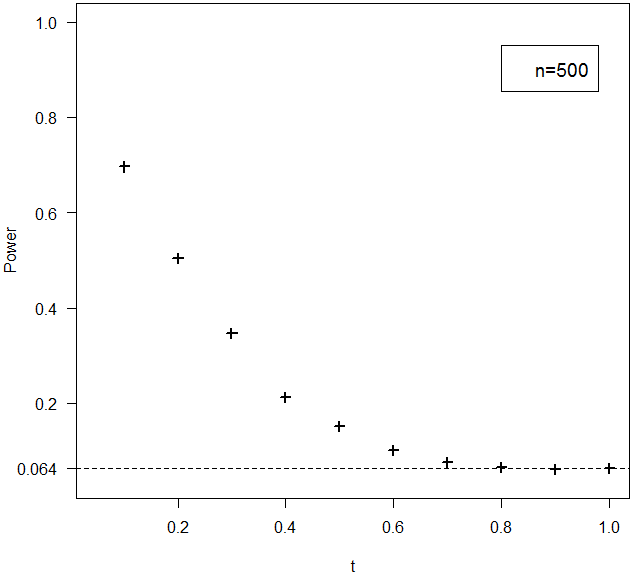} \\
		
			\includegraphics[width=0.46\textwidth,height=4.5cm]{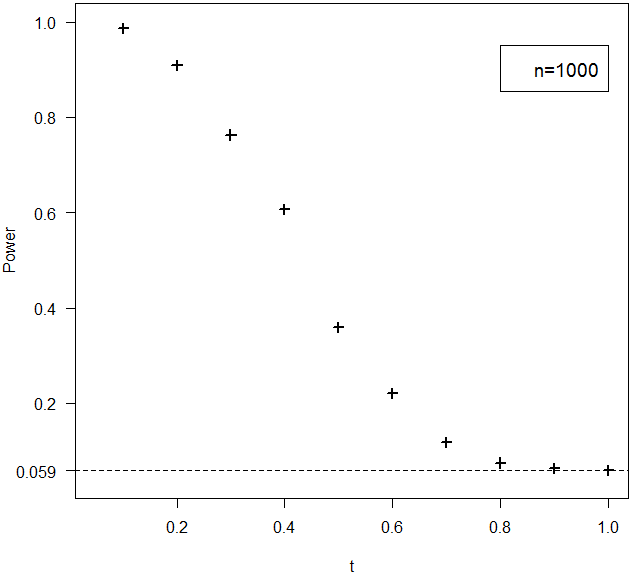} \hfill \includegraphics[width=0.46\textwidth,height=4.5cm]{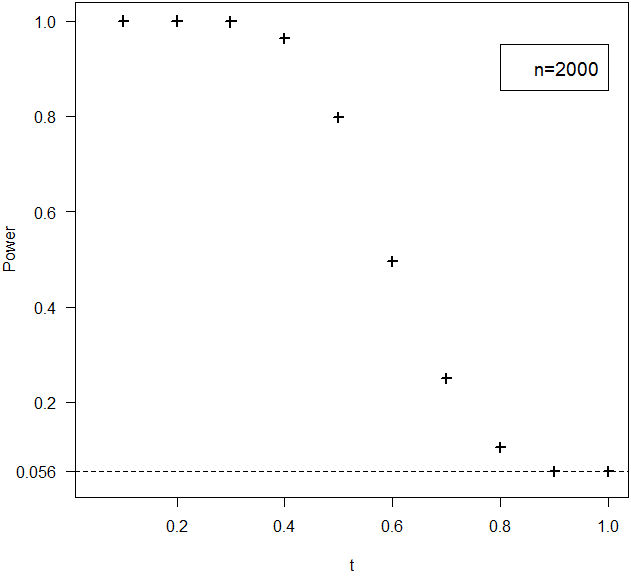}
		\caption{\footnotesize\textit{Test of $\mathcal H_0^{(1)}: \supp(P) \subseteq \supp(P_0)$ in the model $\mathcal M^{(1)} = \{ \vect(M): M \mathbf 1 = \mathbf 1 \}$.} The graphs represent the estimated power of the test under the alternative $P = t P_0 + (1-t) C$, for the four sample sizes $n=200$ (top-left), $n=500$ (top-right), $n=1000$ (bottom-left) and $n=2000$ (bottom-right). The power is estimated for a nominal size $\alpha=0.05$, based on $10^4$ replications.}\label{fig:powerT1}
	\end{figure*}

 To evaluate the power of the test, we generate a Markov chain $X$ with transition matrix $P_t= t P_0 + (1-t)C$ where $C$ is a stochastic matrix with full support drawn randomly beforehand (the same value of $C$ is used during the whole simulation study). We estimate the power for different values of $t$ varying in $(0,1)$ with $0.1$ increments. The empirical reject frequencies under the alternative hypothesis are given in Figure \ref{fig:powerT1}. \\

For small sample sizes, e.g. $n=200$, the estimated type I error probability ($\approx 0.13$ corresponding to $t=1$) is far from the nominal value $\alpha =0.05$. The test is not powerful in this case as we can see that the power remains smaller than $0.4$ for $t \in (0,1)$ and increases in a neighborhood of $1$, which indicates a bias. Nevertheless, the theoretical results are verified for sample sizes as large as $n=500$ and the test is no longer biased. Both the estimated powers and significance level appear to converge to the expected values as $n$ grows.

\subsection*{Test 2.}
We investigate the test procedure for the hypothesis  $$\mathcal H_0^{(2)}: P=P_0,$$ assuming that the support of $P$ is restricted to that of $P_0$. Note that in this case, the first and last row of $P$ are known since they contain only one non-zero element. The model  is given by
$$ \mathcal M^{(2)} = \big\{\vect(M) : \supp(M) \subseteq \supp(P_0) \big\}.  $$
The estimated significance levels are given in Table 2.  \\

\begin{table}[ht]
\begin{center}
\begin{tabular}{|c||c|c|c|c|} \hline 
$ n $ & $200$ & $500$ & $1000$ & $2000$ \\
\hline
\hline
$\hat \alpha$ & $\!\!\! \begin{array}{c} 0.255 \\ (0.0044) \end{array}\!\!\! $ & $ \!\!\! \begin{array}{c} 0.104 \\ (0.0031) \end{array}\!\!\!  $ & $\!\!\! \begin{array}{c} 0.075 \\ (0.0026) \end{array}\!\!\! $ & $\!\!\! \begin{array}{c} 0.057 \\ (0.0023) \end{array}\!\!\! $ \\
\hline
\end{tabular}
\vspace{0.3cm}
\caption{\footnotesize\textit{Test for the hypothesis  $\mathcal H_0^{(2)}: P = P_0$ in the model $\mathcal M^{(2)} = \{ \vect(M): \supp(M) \subseteq \supp(P_0) \}$.} The table contains summary statistics of Monte Carlo simulation results obtained with $P=P_0$ given in Eq. (\ref{ex:p1}), for the four considered sample sizes and for $\mu \sim \mathcal P(1)$. The probability of rejecting $\mathcal H_0^{(1)}$ is estimated for a significance level $\alpha=0.05$, based on $10^4$ replications.}
\end{center}
\label{tab:exemple2}
\end{table}

\begin{figure}[h]
	\centering
		\includegraphics[width=1.09\textwidth,height=4cm]{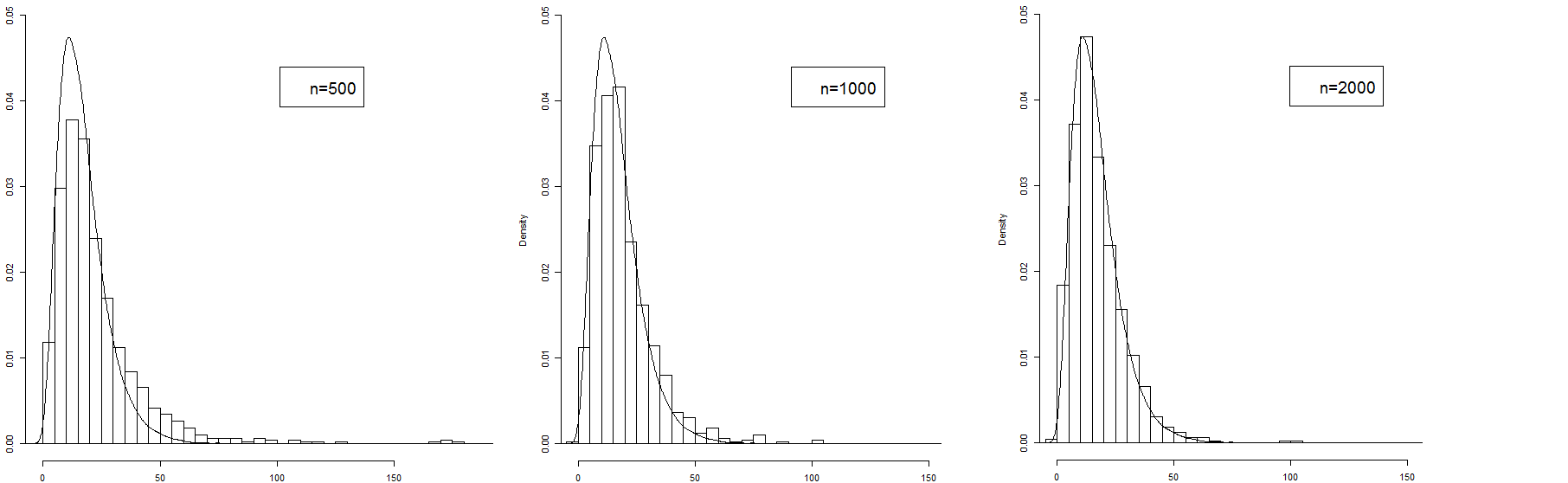}	
\caption{\footnotesize \textit{Histograms of the $10^4$ realizations of the test statistic for the test of $\mathcal H_0^{(2)}: P = P_0$ in the model $\mathcal M^{(2)} = \{ \vect(M): \supp(M) \subseteq \supp(P_0) \}$.} The histograms are obtained for $n=500$ (left), $n=1000$ (center) and $n=2000$ (right) with the density of the limit distribution $\chi^2_G(W)$ added in straight line.}
\label{fig:histoK2}
\end{figure}

The convergence seems slower in this case compared to the previous test procedure, with an estimated significance level not quite in a two standard-deviation range of objective value for $n=2000$. The slow convergence is partly explained by the approximation of the limit distribution by the generalized chi-squared $\chi^2_G(\hat W)$ obtained with the estimated matrix $\hat W$. Nevertheless, the convergence in distribution of the test statistic to the $\chi^2_G(W)$ is clearly observed in Figure \ref{fig:histoK2}.

\subsection*{Test 3.}
We now focus on the distribution  $\mu$ of the random sampling $(\tau_k)_k$
We want to test the hypothesis
$$ \mathcal H_0^{(3)}: \mu = \mu_0 = \mathcal P(1),   $$
in the model $\mathcal M^{(2)}$. This means that we perform the test to see if $\mu$ is a Poisson distribution with parameter $1$, knowing that the support of $P$ is contained in the upper and lower diagonals. This model satisfies the identifiability assumption $\mathbf A1$.

\begin{figure*}[ht]

	\includegraphics[width=0.46\textwidth,height=4.5cm]{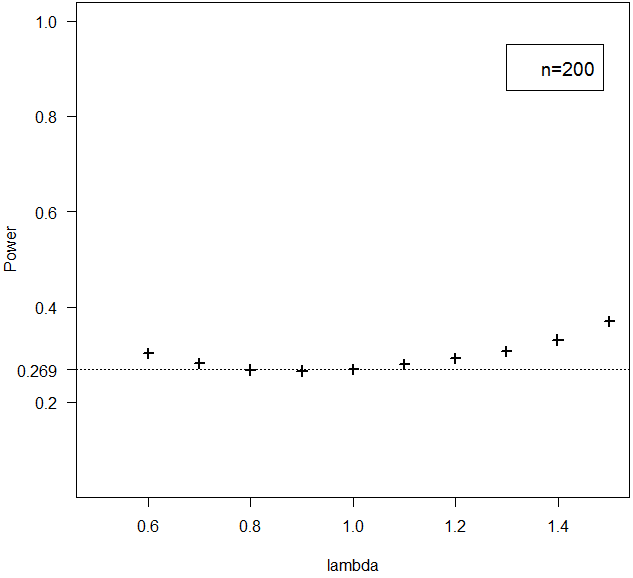}  \hfill \includegraphics[width=0.46\textwidth,height=4.5cm]{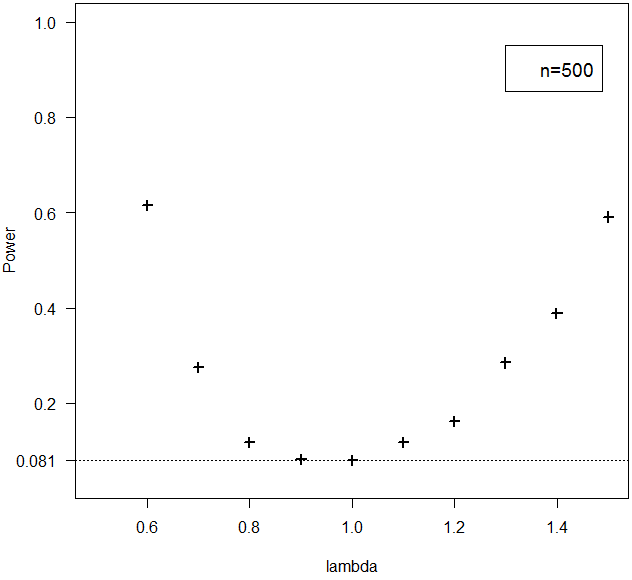} \\
	
			\includegraphics[width=0.46\textwidth,height=4.5cm]{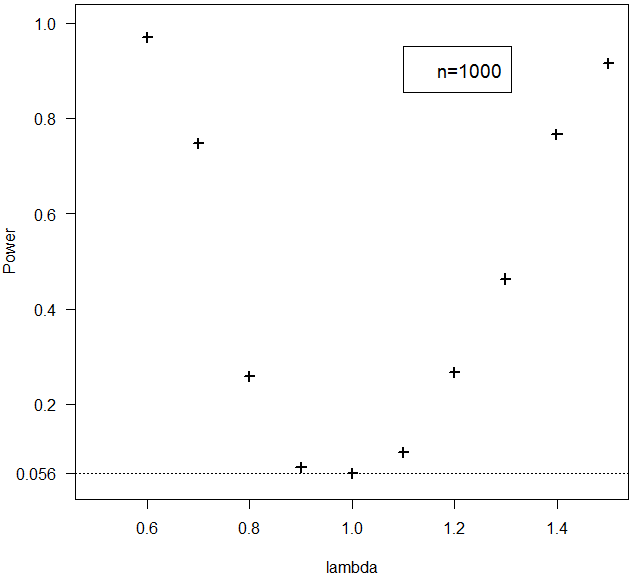} \hfill \includegraphics[width=0.46\textwidth,height=4.5cm]{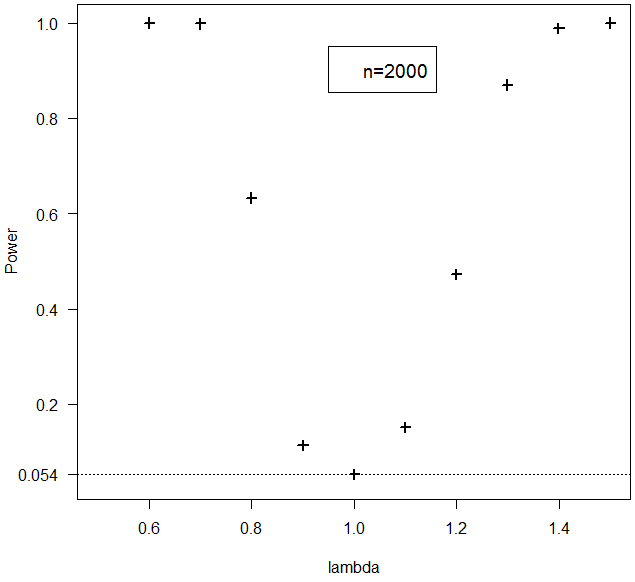}
			
\caption{\footnotesize\textit{Test for the hypothesis $\mathcal H_0^{(3)}: \mu = \mathcal P(1)$ in the model $\mathcal M^{(2)} = \{ \vect(M): \supp(M) \subseteq \supp(P_0) \}$.} The graphs represent the estimated power under the alternative $\mu = \mathcal P(\lambda)$, for $\lambda$ in $(0.5,1.5)$ and for $n=200$ (top-left), $n=500$ (top-right), $n=1000$ (bottom-left) and $n=2000$ (bottom-right). The power is estimated for a nominal size $\alpha=0.05$, based on $10^4$ replications.}
\label{tab:fig2}
\end{figure*}

We estimate  the probability of rejecting the null hypothesis when the time gaps $\tau_k$ are distributed from a Poisson distribution $\mathcal P(\lambda) $, with the parameter $\lambda$ varying in $(0.5,1.5)$ with $0.1$ increments. Table \ref{tab:exemple3} provides the empirical rejection frequency under the null hypothesis $\lambda=1$. Figure \ref{tab:fig2} gives the  estimated power of the test as function of the parameter $\lambda$.

\begin{table}[h]
\begin{tabular}{|c||c|c|c|c|} \hline 
$ n $ & $200$ & $500$ & $1000$ & $2000$ \\
\hline
\hline
$\hat \alpha$ & $\!\!\! \begin{array}{c} 0.269 \\ (0.0044) \end{array}\!\!\! $ & $ \!\!\! \begin{array}{c} 0.081 \\ (0.0027) \end{array}\!\!\!  $ & $\!\!\! \begin{array}{c} 0.056 \\ (0.0023) \end{array}\!\!\! $ & $\!\!\! \begin{array}{c} 0.054 \\ (0.0023) \end{array}\!\!\! $ \\
\hline
\end{tabular}
\vspace{0.3cm}
\caption{\footnotesize\textit{Test for the hypothesis  $\mathcal H_0^{(3)}: \mu = \mathcal P(1)$ in the model $\mathcal M^{(2)} = \{ \vect(M): \supp(M) \subseteq \supp(P_0) \}$.} The table contains summary statistics of Monte Carlo simulation results obtained with $P=P_0$ for the four sample sizes under the null hypothesis. The nominal size of the test is $\alpha=0.05$, based on $10^4$ replications.}
\label{tab:exemple3}
\end{table}

\vspace{-0.3cm}

The power of the test is not satisfactory for $n=200$ but it increases significantly with $n$. Although the time gaps are not observed, this simulation confirms that performing a goodness-of-fit test on their distribution is still possible in this framework. In particular, the intensity of the process assuming a Poisson distribution can be recovered from the indirect observations. The hypothesis test turns out to be rather satisfactory provided the number of observations is sufficient.

\section{Conclusion}\label{sec:concl}
We develop an hypothesis testing methodology on the transition matrix $P$ and the distribution $\mu$ of time gaps in a Markovian model with random sampling. The original contribution consists in the construction of a testing procedure for affine hypothesis on $P$ in this framework as well as a test on the distribution $\mu$ of the observation times. This setting requires additional information on $P$ to make the model identifiable. The information is expressed in the form of affine constraints on $P$, thus extending the framework considered in \cite{fytp14} in which $P$ is assumed to be sparse. \\

We show that the different forms of information on $P$, formalized via affine constraints, can be sufficient to recover the true distribution $\mu$ of the time gaps, or at least to build a consistent testing procedure on specific values of $\mu$. The simulation study confirms that a goodness-of-fit test on $\mu$ can be carried out successfully despite the fact that the time gaps $\tau_k$ are not observed.

\section{Appendix}\label{sec:app}

\noindent \textit{Proof of Theorem \ref{th1}.} We use the fact that, for $\mathcal A$ an affine subset of an Euclidean space,
$$ \inf_{a \in \mathcal A} \Vert a \Vert^2 = \Vert (\id - \Pi_{\mathcal A_{\lin}}) a_0 \Vert^2,   $$
for all $a_0 \in \mathcal A$, with $\Pi_{\mathcal A_{\lin}}$ the orthogonal projector onto $\mathcal A_{\lin} = \mathcal A - \mathcal A$. Define the linear sets $\hat E = \im(\Delta(\hat Q)\Phi)$ and $\hat E_0 = \im(\Delta(\hat Q)\Phi_0)$,
we get
$$  \inf_{m \in \mathcal M_0} \Vert \Delta(\hat Q) m \Vert^2 - \inf_{m \in \mathcal M} \Vert \Delta(\hat Q) m \Vert^2 = \Vert (\id - \Pi_{\hat E_0}) \Delta(\hat Q) m_0 \Vert^2 - \Vert (\id - \Pi_{\hat E}) \Delta(\hat Q) m_0 \Vert^2 $$
for any $m_0$ in $\mathcal M_0$ and in particular, for $m_0 = p$ under $\mathcal H_0$. Since $\hat E_0 \subseteq \hat E$, the Pythagorean theorem yields, under $\mathcal H_0$,
$$ \inf_{m \in \mathcal M_0} \Vert \Delta(\hat Q) m \Vert^2 - \inf_{m \in \mathcal M} \Vert \Delta(\hat Q) m \Vert^2 =  \Vert (\Pi_{\hat E} - \Pi_{\hat E_0} )\Delta(\hat Q) p \Vert^2 =  \Vert \Pi_{\hat F}\Delta(\hat Q) p \Vert^2,  $$
setting $\hat F = \hat E \cap \hat E_0^\perp$. Write $\Pi_{\hat F } \Delta(\hat Q) p  = \Pi_{F}\Delta(\hat Q) p  + (\Pi_{\hat F} - \Pi_F) \Delta(\hat Q) p $, we have
$$ \sqrt n \ \Pi_F \Delta(\hat Q) p = - \sqrt n \ \Pi_F \Delta(P) \hat q \overset{d}{\longrightarrow} \mathcal N \big(0,\Pi_F \Delta(P) \Sigma \Delta(P)^\top \Pi_F \big),   $$
yielding $ n \Vert \Pi_F \Delta(\hat Q) p \Vert^2  \overset{d}{\longrightarrow} \chi^2(W)$. So, to prove the result, it remains to show that 
$$n \big( \Vert \Pi_F \Delta(\hat Q) p \Vert^2 - \Vert \Pi_{\hat F} \Delta(\hat Q) p \Vert^2 \big) = o_P(1).$$
Since $\sqrt n \Vert \Delta(\hat Q) p \Vert = O_P(1)$, it suffices to show that $\Pi_{\hat F}$ converges in probability to $\Pi_F$ in view of 
\begin{eqnarray*} n \big( \Vert \Pi_F \Delta(\hat Q) p \Vert^2 - \Vert \Pi_{\hat F} \Delta(\hat Q) p \Vert^2 \big) & \leq & n \Vert (\Pi_F - \Pi_{\hat F}) \Delta(\hat Q) p \Vert \ \big( \Vert \Pi_F \Delta(\hat Q) p \Vert + \Vert \Pi_{\hat F} \Delta(\hat Q) p \Vert \big) \\
& \leq & \Vvert \Pi_F - \Pi_{\hat F} \Vvert \times 2 n \Vert \Delta(\hat Q) p \Vert^2,
\end{eqnarray*}
where $\Vvert . \Vvert$ denotes the operatorn norm. Recall that $\Pi_{\hat F} = \Pi_{\hat E} - \Pi_{\hat E_0}$ with $\hat E = \im(\Delta(\hat Q) \Phi)$ and $\hat E_0 = \im(\Delta(\hat Q) \Phi_0)$. By \textbf{A1}, $\Delta(Q) \Phi_0$ is of full rank and
$$\Pi_{\hat E_0} = \Delta(\hat Q) \Phi_0 \left(\Phi_0^\top \Delta(\hat Q)^\top \Delta(\hat Q) \Phi_0\right) ^{-1}\Phi_0^\top \Delta(\hat Q)^\top $$ 
clearly converges in probability to $\Pi_{E_0}$. Moreover, we know by \textbf{A2} that $\dim(\hat E) = \im(\Delta(\hat Q) \Phi) \leq\dim(E)$, so that the convergence of $\Delta(\hat Q) \Phi$ to $\Delta(Q) \Phi$ is sufficient for $\Pi_{\hat E}$ to converge to $\Pi_E$. Thus, $\Pi_{\hat F}$ converges in probability to $\Pi_F$ as $n \to \infty$ which ends the proof. \\

\vskip 1cm 

\noindent \textit{Proof of Corollary \ref{coro}.} By continuity of the map $(\epsilon,W) \mapsto \epsilon^\top W \epsilon$, Slutsky's lemma ensures the convergence in distribution of $\chi^2_G(\hat W)$ to $\chi^2_G(W)$ and therefore, the convergence of the estimated quantile $\hat u_{1-\alpha}$ to the true quantile $u_{1-\alpha}$ of the $\chi^2_G(W)$ distribution. It follows that $\lim_{n \to \infty} \mathbb P(S > \hat u_{1-\alpha}) = \mathbb P(S > u_{1-\alpha}) = \alpha$. Now, if $\forall p' \in \mathcal M_0: \ \Delta(Q) p' \neq 0$, we show easily that $ \inf_{m \in \mathcal M_0} \Vert \Delta(\hat Q) m \Vert^2$ does not converge to zero as $n$ grows to infinity, while we still have
$$ \inf_{m \in \mathcal M} \Vert \Delta(\hat Q) m \Vert^2 \leq \Vert \Delta(\hat Q) p \Vert^2 \overset{\mathbb P}{\longrightarrow} 0. $$
Hence, $ S = n \big( \inf_{m \in \mathcal M_0} \Vert \Delta(\hat Q) m \Vert^2 - \inf_{m \in \mathcal M} \Vert \Delta(\hat Q) m \Vert^2 \big)$ diverges in probability in this case, and thus the result follows. 

\vskip 1cm 

\noindent \textit{Proof of Proposition \ref{test}.} To show that the matrix $\Gamma$ exists, it suffices to show that the series is normally convergent. We use that $\Vvert (P^{j-1})^\top \otimes P^{k-j} \Vvert \leq \Vvert P^{j-1} \Vvert \ \Vvert P^{k-j} \Vvert \leq 1$ to get
$$  \sum_{k \geq 1} \bigg( \sum_{j=1}^{k} \Vvert (P^{j-1})^\top \otimes P^{k-j} \Vvert \bigg) \mu(k) \leq \sum_{k \geq 1} k \mu(k) = \mathbb E_\mu(\tau),    $$
which is finite by assumption. We now compute the differential of $g_\mu$ at $P$. For $H \in \mathbb R^{N \times N}$ such that $\Vvert H \Vvert \leq 1$, we have
\begin{eqnarray*} \lim_{t \to 0} \frac 1 t \big( G_\mu(P + t H) - G_\mu(P)  \big) & = & \lim_{t \to 0} \ \frac 1 t \ \sum_{k \geq 1} \big[ (P + t H)^k - P^k \big] \mu(k)\\
& = & \lim_{t \to 0} \frac 1 t \ \bigg[  \sum_{k \geq 1} t \bigg( \sum_{j=1}^k P^{j-1} H P^{k-j}   \bigg) \mu(k) + o(t) \bigg] \\
& = & \sum_{k \geq 1} \bigg( \sum_{j=1}^k P^{k-j} H P^{j-1} \bigg) \mu(k).
\end{eqnarray*}
Applying the vectorization yields, in view of $\vect(ABC) = (C^\top \otimes A) \vect(B)$,
\begin{eqnarray*} \lim_{t \to 0} \frac 1 t \big( g_\mu(P + t H) - g_\mu(P)  \big) & = & \sum_{k \geq 1} \bigg( \sum_{j=1}^k \vect \big(P^{k-j} H P^{j-1} \big)   \bigg) \mu(k) \\
 & = & \bigg[\sum_{k \geq 1} \bigg( \sum_{j=1}^{k} (P^{j-1})^\top \otimes P^{k-j}\bigg) \mu(k) \bigg] \vect(H) \\
 & = & \Gamma \vect(H).
\end{eqnarray*}
It follows from Cramer's theorem that, 
$$ \sqrt n \big( g_\mu(\hat P) - q \big) = \sqrt n \ \Gamma \big( \hat p - p \big) + o_P (1).  $$
Using \eqref{eqhatp}, we deduce that $ \sqrt n ( g_\mu(\hat P) - q ) =  \sqrt n \ \Gamma B ( \hat q - q ) + o_P ( 1)$, which leads to
$$ \sqrt n \ \big( \hat q - g_\mu(\hat P) \big) = \sqrt n \ \big( \id - \Gamma B \big) (\hat q - q) + o_P (1),  $$
and the result follows.\\

\bibliographystyle{acm}
\bibliography{refs}

\end{document}